\documentclass[11pt]{article}

\usepackage{latexsym,verbatim,ifthen,graphicx,mathrsfs,tabularx}
\usepackage[english]{babel}
\usepackage[utf8]{inputenc} 	
\usepackage{amssymb}
\usepackage{amsmath}
\usepackage{amsthm}	
\usepackage{tikz-cd}
\usepackage{mathtools}
\usepackage{microtype}
\usepackage[overload]{empheq}
\usepackage{todonotes}
\usepackage{macros}
\usepackage{calc}
\usepackage{comment}
\newlength\myheight
\newlength\mydepth
\settototalheight\myheight{Xygp}
\usetikzlibrary{matrix,arrows,decorations.pathmorphing}
\usepackage[all]{xy}				
\usepackage{hyperref}
\hypersetup{
	colorlinks=true,
	allcolors=blue,
}   		
\usepackage{anysize}
\usepackage{cancel}									
\marginsize{2.5cm}{2.5cm}{2cm}{2cm}
\usepackage{enumerate}		

 \usepackage{booktabs}

\usepackage{mathrsfs}									
\usepackage{fourier}			
\usepackage[Sonny]{fncychap}	
\usepackage{fancyhdr}

\newtheorem{theorem}{Theorem}[section]

\newtheorem{proposition}[theorem]{Proposition}
\newtheorem{corollary}[theorem]{Corollary}
\theoremstyle{definition}
\newtheorem{example}[theorem]{Example}

\theoremstyle{definition}
\newtheorem{remark}[theorem]{Remark}

\newtheorem{definition}[theorem]{Definition}

\begin{document}

\title{Determinant Massey products}

\author{Oisín Flynn-Connolly, José M. Moreno-Fernández, Fernando Muro}

\date{}
\maketitle
\abstract{We introduce the determinant Massey product, 
a secondary operation on the cohomology of commutative dg algebras which is 
reminiscent of the classical Massey product.
We prove some of the elementary properties of this product, compute several examples,
and relate them to classical and matric Massey products.
When 2 is invertible, a determinant Massey product is a matric Massey product,
but over $\mathbb F_2$, it is not. In this case, we produce 
a finite-dimensional cdga whose Harrison cohomology 
does not inject into its Hochschild cohomology, 
and a divided powers algebra which is formal as an 
associative dg algebra, but not as a commutative one,
with an explicit determinant Massey product witnessing the failure.
Over $\mathbb Q$, we use determinant Massey products to construct a 
simply-connected coformal space with Koszul cohomology ring and vanishing ordinary 
Massey products of all orders, which is not formal.}

\section{Introduction}

Ever since W. S. Massey introduced what nowadays are called \emph{Massey products} \cite{UeharaMassey57,Massey58}, 
a secondary operation on the (co)homology of differential graded associative algebras,
these have come to play a major role in algebraic topology and spread to many other seemingly unrelated areas.
The existence of these and related higher-order products is due to the vanishing of certain equations
that follow from associativity, commutativity, and more complex operadic operations \cite{muro23, operadicmassey}. 
Their first application was the Jacobi identity for Whitehead products \cite{UeharaMassey57},
and they detect that the Borromean rings are non-trivially linked \cite{Mas69,Mas98}.
Subsequent applications to algebraic topology abound. 
For example, Longoni and Salvatore used them to show 
that configuration spaces do not reflect homotopy equivalences \cite{Longoni05},
and Livernet to show that the Swiss cheese operad is not formal \cite{livernet15}. 
They have also been extensively studied in the context of the Adams spectral sequence 
and the computation of the stable homotopy groups of spheres \cite{Isaksen20, Kochman80}. 
Recently, their behaviour in the Galois cohomology of number fields has become of major interest in arithmetic geometry
\cite{EM17, HW15, GMT18, HW19, MT15b, MT16}, including playing a role in constructions of integral points on elliptic curves \cite{Minhyong10}. 
They have appeared as far afield as in the analysis of differential equations \cite{Reyesetal} and mathematical physics \cite{Bruinsma}.
Despite the recent interest in $\infty$-algebras,
which are marvellous theoretical tools,
they are not exactly the same thing as Massey products \cite{jose20, operadicmassey}, 
which are better suited to practical computations and the construction of counterexamples.

\smallskip

A significant limitation of classical Massey products is that they rely on initial monomial vanishing conditions in cohomology, 
such as \( ab = 0 \). To address this, May introduced a matrix-based formalism to handle vanishing conditions given by polynomials,
expressed as \( \sum_i a_i b_i = 0 \), 
leading to the development of \emph{matric Massey products} \cite{May69}.
These also play a major role in the most recent computations of stable homotopy groups of spheres \cite{gheorge21}. 
In this paper, we introduce \emph{determinant Massey products}.
As the name suggests, rather than matrix multiplication, 
these are based on taking determinants of matrices.
They are defined on commutative dg algebras, 
where the determinant of a matrix with two equal columns vanishes.
When $2$ is invertible, they are matric Massey products of a special shape
(Proposition~\ref{prop: connection to matric Massey products}), 
but over $\mathbb F_2$ they are not. 
There, they detect commutative structure that is invisible to associative invariants,
see Section~\ref{sec: ass vs com in char 2}.

\smallskip

In \cite{campos24}, Campos, Petersen, Robert-Nicoud and Wierstra show 
that two commutative dg algebras over a field of characteristic zero that 
are quasi-isomorphic as associative dg algebras are also quasi-isomorphic as commutative dg algebras. 
The first author showed, using Frobenius-type obstructions, that this fails in positive characteristic \cite{Oisin--ObstructionTheory}. 
Determinant Massey products cannot detect this failure when $2$ is invertible, 
as matric Massey products are invariants of the associative homotopy type. 
Over $\mathbb F_2$ they can detect this failure at the level of deformation complexes. We give examples of a graded commutative algebra whose Harrison cohomology does not 
inject into its Hochschild cohomology (Proposition \ref{prop: Harrison not injective}),
and a divided powers algebra that is formal as an associative dg algebra but not as a commutative one (Proposition \ref{prop: dividedpowers}), significantly strengthening the results of \cite{Oisin--ObstructionTheory}.
Over $\mathbb Q$, the simplest possible algebra with a non-vanishing determinant Massey product is also intrinsically interesting. 
It provides an example of a simply-connected coformal space with Koszul cohomology ring 
and vanishing ordinary Massey products of all orders, that is not formal (Proposition~\ref{prop: berglund space})

\smallskip 

\subsubsection*{Structure of the article}

After this introduction, we establish our notational conventions and the specific definition of determinants
for matrices over commutative graded algebras.
In Section \ref{sec: determinant Massey products},
we introduce the determinant triple Massey product, prove its fundamental properties, 
and study its behaviour under morphisms of differential graded algebras. 
In Section~\ref{sec: connections to classical and matric Massey products}, 
we relate it to classical and matric Massey products, 
and in Section~\ref{sec: applications} we collect applications and examples.
Finally, in Section \ref{sec: Higher-order determinant Massey products},
we sketch how to generalise our construction to define higher-arity determinant Massey products.

\subsubsection*{Notation and conventions}

We use the cohomological convention for differential graded algebras, i.e., differentials raise degree by one.
If $f:A\to B$ is a morphism of dg algebras, we consistently denote by $f^*$ the map induced by $f$ in cohomology.
The class of a cocycle $x$ will be denoted by $\overline{x}$, and its degree by $|x|$.
All the algebras in this paper are taken over a commutative unital ground ring.
We specify its characteristic whenever it matters. We abbreviate commutative dg algebra to cdga.
Although part of what is done in this paper admits a generalization to associative algebras,
we will assume all the algebras to be graded-commutative unless otherwise specified.

\subsubsection*{On determinants}

Here, we lay down our convention for the determinant of square matrices 
over commutative graded algebras.  
We mainly aim for a closed formula to help us keep track of the combinatorics involved in our determinant Massey products.

\smallskip

Let $A$ be a commutative graded algebra. 
The determinant of a square matrix $M$ with $A$-coefficients is defined as follows.
If $M=(a)$ is a size 1 matrix, its determinant is $a$.
For size 2 matrices, 
we use the Koszul sign graded determinant: 
\begin{equation*}
\begin{vmatrix}
a & b \\
c & d 
\end{vmatrix} = ad - (-1)^{|b||c|} cb = ad - bc.
\end{equation*}
For matrices of size $n\geq 3$, we inductively expand along the first column by the formula
\begin{equation*}
|M|= \sum_{i=1}^{n} (-1)^{i +1}\cdot \varepsilon_i\cdot \ a_{i1} \cdot \ |M^{(i,1)}|.
\end{equation*}
Here, writing $a_i=a_{i1}$ for the entries of the first column,
$\varepsilon_i$ is the Koszul sign associated to the permutation that moves $a_i$ to 
the first place:
\begin{equation*}
\left(a_1,a_2,\ldots,a_i,\ldots,a_n\right) 
\mapsto (-1)^{|a_i|\left(\sum_{j=1}^{i-1}|a_j|\right)}\left(a_i,a_1,a_2,\ldots,\widehat{a_i},\ldots,a_n\right), 
\end{equation*}
and $M^{(i,1)}$ is the square submatrix of $M$ obtained by removing its $i$-th row and first column.
Note that the formula for size 2 is not the case $n=2$ of this expansion, 
which would give $ad-(-1)^{|a||c|}cb$. The two agree whenever $a$ and $b$ have the same parity,
but not in general, and Remark~\ref{rmk: parity} uses the formula for size 2 in a case where they do not.
For example, the explicit formula for size 3 matrices is the following: 
\begin{equation*}
\begin{vmatrix}
a & d & g \\
b & e & h \\
c & f & i 
\end{vmatrix} =  a \left( e i - hf \right) - (-1)^{|b||a|} b \left( d i - gf \right)
+ (-1)^{|c||a| + |c||b|} c \left( d h - ge \right).
\end{equation*}
If $A$ is concentrated in even degrees, then the signs in the determinant expansion become the usual ones.

\section{Determinant Massey products}
\label{sec: determinant Massey products}

In this section, we introduce the determinant (triple) Massey products.
After proving some of its basic properties, we give concrete applications and examples.
Since triple determinant Massey products are potentially the most useful in practice,
we begin the paper with a detailed study of them, and address the higher-order case separately in the last section of the paper.

\smallskip

We begin with the definition of the determinant Massey products.

\begin{definition}
\label{def: determinant Massey product}
Let $\bar a,\bar b,\bar c,\bar x,\bar y, \bar z\in H^*(A)$ be homogeneous classes in the cohomology of a commutative 
dg algebra $A=\left(A,d\right)$.
A \emph{defining system} for the determinant Massey product of
$\bar a, \bar b, \bar c$ over $\bar x,\bar y,\bar z$
is a choice of elements 
\begin{equation*}
\mD = \left\{a,b,c,\ \ x,y,z,\ \ \zeta_a,\zeta_b,\zeta_c\right\}
\end{equation*}
such that $a,b,c,x,y,z$ are cocycle representatives of $\bar a, \bar b, \bar c, \bar x,\bar y, \bar z$, respectively,
and $\zeta_a,\zeta_b,\zeta_c$ are homogeneous elements satisfying
\begin{equation*}
    \label{ecu: bounding assumptions}
    d\zeta_a=yc-bz, \quad \quad 
d\zeta_b=xc-az, \quad \quad 
d\zeta_c=xb-ay.
\end{equation*}
Furthermore, we assume that
\footnote{The first two equalities below imply the third, so in practice,
it suffices to check those.}
\begin{equation*}
|y|+|c| = |z|+|b|, \qquad |x|+|c| = |z| + |a|,
\qquad \textrm{and} \qquad |x|+|b|=|y|+|a|,
\end{equation*}
equivalently, the three differences $|x|-|a| \;=\; |y|-|b| \;=\; |z|-|c|$ coincide,
and the following congruence relation mod 2 is satisfied:
\begin{equation*}
|x|-|a| \;=\; |y|-|b| \;=\; |z|-|c| \equiv 0 \pmod 2.
\end{equation*}
This is to say that $x$ and $a$ (respectively, $y$ and $b$, and $z$ and $c$) have the same parity.
In this case, the element 
\begin{equation}
    \label{ecu: determinant Massey product cocycle}
    U_{\mD} \,=
    (-1)^{|a|}a\zeta_a - (-1)^{|b| + |b||a|}b \zeta_b + (-1)^{|c|+ |c||a| + |c||b| }c\zeta_c
\end{equation}
is a cocycle whose class $\overline{U_{\mD}}$
is the \emph{determinant Massey product of
$\bar a, \bar b, \bar c$ over $\bar x,\bar y,\bar z$} associated to the defining system $\mD$.
We denote by $\langle \bar a,\bar b,\bar c \ :\ \bar x, \bar y, \bar z\rangle_{\det}$ the (possibly empty) set of cohomology classes
formed by all possible choices of defining systems, 
\begin{equation*}
\langle \bar a,\bar b,\bar c \ :\ \bar x, \bar y, \bar z\rangle_{\det} = \left\{ \overline{U_{\mD}} \ \mid \
\mD = \left\{a,b,c,\ \ x,y,z,\ \ \zeta_a,\zeta_b,\zeta_c\right\} \textrm{ is a def. system} \right\}
\subseteq H^{|a|+|y|+|c|-1}\left(A\right).
\end{equation*}
The determinant Massey product is \emph{defined} if the set $\langle \bar a,\bar b,\bar c \ :\ \bar x, \bar y, \bar z\rangle_{\det}$ is non-empty,
that is, if there is some defining system as above.
\end{definition}

Several observations are in order.
First, $U_{\mD}$ is a cocycle:
\begin{equation*}
dU_{\mD}=a(yc-bz)-(-1)^{|a||b|}b(xc-az) +(-1)^{|a||c|+|b||c|}c(xb-ay). 
\end{equation*} 
The two terms involving $abz$ cancel by graded commutativity. 
The remaining ones cancel in pairs, because $|y|\equiv|b|$ and $|x|\equiv|a|$.
Below, we interpret this computation in terms of determinants.

Second, when mentioning determinant Massey products, 
we will omit the part "of $\bar a,\bar b,\bar c$ over $\bar x,\bar y,\bar z$"
and take it as implicit. 
However, it is important that $\bar a,\bar b,\bar c$ comes before $\bar x, \bar y, \bar z$, that is, the definition is not 
completely symmetric on $\bar a,\bar b,\bar c$ with respect to $\bar x, \bar y, \bar z$; this will be explained in 
Remark \ref{rmk: Determinant Massey products appear in pairs}.

\smallskip

Our first observation is that the determinant Massey product does not 
depend on the choice of cocycle representatives.

\begin{proposition}
\label{prop: determinant massey doesnt depend on cocycle representatives}
The determinant Massey product $\langle \bar a,\bar b,\bar c \ :\ \bar x, \bar y, \bar z\rangle_{\det}$ 
does not depend on the choice of cocycle representatives.
That is, if $a,b,c,x,y,z$ represent $\bar a,\bar b,\bar c,\bar x,\bar y,\bar z$, 
then every class in 
$\langle \bar a,\bar b,\bar c \ :\ \bar x, \bar y, \bar z\rangle_{\det}$ is
$\overline{U_{\mD}}$ for a defining system $\mD$ whose representatives are $a,b,c,x,y,z$. 
\end{proposition}

\begin{proof}
Let $\overline{U_{\mD}}$ be a determinant Massey product 
with defining system $\mD=\left\{a,b,c,\ \ x,y,z,\ \ \zeta_a,\zeta_b,\zeta_c\right\}$.
The elements $\bar a,\bar b,\bar c$ play a 
slightly different role than $\bar x,\bar y,\bar z$ in this operation, so the proof breaks in two cases. 
It suffices to replace one representative at a time. 
We give the details when replacing $a$ by $a'=a+dl$, and $x$ by $x'=x+de$.
The other cases are analogous.
\begin{enumerate}
    \item[{(i)}] Assume we replace $a$ by $a' = a +dl$. 
    Let us prove that the set $\{a',b,c\}$ extends to a defining system
$\mD'=\left\{a',b,c,\ \ x,y,z,\ \ \zeta_a,\zeta_b',\zeta_c'\right\}$ 
such that the associated cocycle $U_{\mD'}$ is cohomologous to $U_\mD$.
Indeed, it suffices to take
\begin{equation*}
\zeta_b' = \zeta_b -lz \qquad \textrm{and} \qquad \zeta_c' = \zeta_c - ly,
\end{equation*}
and leave the rest of elements as in the original defining system.
In this case, recalling that $|l|=|a|-1$ and $|y|\equiv|b|$, we find that the classes are cohomologous:
\begin{align*}
    U_{\mD'}  &= (-1)^{|a|}a'\zeta_a - (-1)^{|b| + |b||a|}b \zeta_b' + (-1)^{|c|+ |c||a| + |c||b| }c\zeta_c' \\[0.25cm]
        &= (-1)^{|a|}(a+dl)\zeta_a - (-1)^{|b| + |b||a|}b \left(\zeta_b -lz\right) 
        + (-1)^{|c|+ |c||a| + |c||b| }c\left(\zeta_c - ly\right) \\[0.25cm]
        &= U_\mD + (-1)^{|a|} \left(dl\right) \zeta_a + (-1)^{|b|+|b||a|} blz - (-1)^{|c|+|c||a| + |c||b|} cly  \\[0.25cm]
        &= U_\mD + (-1)^{|a|} \left(dl\right) \zeta_a +lbz - lyc \quad =\quad   U_\mD + d\left( (-1)^{|a|}l\zeta_a \right). 
\end{align*}
\item[{(ii)}] Assume we replace $x$ by $x' = x +de$. 
We will modify the defining system in such a way that we obtain exactly the same cocycle.
We have:
\begin{align*}
    d\zeta_b    &= xc-az 
    = \left(x'-de\right)c-az 
    = x'c-az - (de)c.
\end{align*}
Since $c$ is a cycle, the element $\zeta_b' = \zeta_b + ec$ is such that  
\begin{equation*}
d\zeta_b' = x'c-az.
\end{equation*}
Similarly, the element $\zeta_c' = \zeta_c + eb$ is such that 
$d\zeta_c' = x'b-ay$.
This yields the determinant Massey product representative
\begin{equation*}
U' \,=(-1)^{|a|}a\zeta_a - (-1)^{|b| + |b||a|}b \zeta_b' + (-1)^{|c|+ |c||a| + |c||b| }c\zeta_c'.
\end{equation*}
It follows from commutativity that the difference between the original cocycle and the one above is zero.
\end{enumerate}
\end{proof}

Moving a level up in the choice of elements of a defining system does not behave as well: 
it produces indeterminacies.
Before dealing with these, we explain the connection to determinants.

\smallskip

Let $\mD$ be a defining system producing the cocycle 
$$U_{\mD} \,=(-1)^{|a|}a\zeta_a - (-1)^{|b| + |b||a|}b \zeta_b + (-1)^{|c|+ |c||a| + |c||b| }c\zeta_c.$$
Then the boundary of $U_\mD$
is formally the determinant of the following matrix when recursively expanding along the first column: 
$$ \begin{pmatrix}
a & x & a \\
b & y & b \\
c & z & c 
\end{pmatrix}.$$
That is, 
\begin{align*}
dU_{\mD} &=a\left(d\zeta_a\right) - (-1)^{|a||b|}b(d\zeta_b) + (-1)^{|a||c|+|b||c|}c(d\zeta_c)  \\[0.25cm]  
&= a\left(yc-bz\right) -(-1)^{|b||a|}b\left(xc-az\right)
+(-1)^{|c||a|+|c||b|}c\left(xb-ay\right)\\[0.25cm] 
&= \ 
a\underbrace{\begin{vmatrix}
y & b \\
z & c 
\end{vmatrix}}_{d\zeta_a}
\ -\ (-1)^{|a||b|}\, b\underbrace{\begin{vmatrix}
x & a \\
z & c 
\end{vmatrix}}_{d\zeta_b}
\ +\ (-1)^{|a||c|+|b||c|}\, c \underbrace{\begin{vmatrix}
x & a \\
y & b 
\end{vmatrix}}_{d\zeta_c} 
\ =\ \begin{vmatrix}
a & x & a \\
b & y & b \\
c & z & c 
\end{vmatrix}.
\end{align*}
Remarkably, the elements $\zeta_a,\zeta_b,$ and $\zeta_c$ bound
exactly the minors occurring in this expansion:
\begin{equation*}
d\zeta_a = \begin{vmatrix}
y & b \\
z & c 
\end{vmatrix}, \qquad \quad 
d\zeta_b = 
\begin{vmatrix}
x & a \\
z & c 
\end{vmatrix}, \qquad \quad 
d\zeta_c =\begin{vmatrix}
x & a \\
y & b 
\end{vmatrix}.
\end{equation*}

The vanishing of these minors in cohomology
gives precisely the conditions imposed to the determinant Massey product to exist,
and it is the reason for their name.

\begin{remark}
\label{rmk: order two}
The construction above also makes sense for two classes, where it gives a primary operation.
This is computed explicitly by letting $\bar a,\bar b,\bar x,\bar y\in H^*(A)$ be homogeneous classes with $|x|-|a|=|y|-|b|$.
We can then define an operation by taking the determinant of the $2\times 2$ matrix
with columns $(\bar a,\bar b)^T$ and $(\bar x,\bar y)^T$,
\begin{equation*}
\langle \bar a,\bar b : \bar x,\bar y\rangle_{\det}
:=\left\{\begin{vmatrix} \bar a & \bar x\\ \bar b & \bar y\end{vmatrix}\right\}
=\left\{\bar a\,\bar y-(-1)^{|b||x|}\,\bar b\,\bar x\right\}.
\end{equation*}
It is obvious that this product is always defined, consists of a single class,
and is preserved by morphisms of cdgas.
It is to the determinant Massey product as the cup product $\bar a\,\bar b$ is to classical Massey products.
Namely, up to sign, the minors above represent the order-two products
$\langle \bar b,\bar c : \bar y,\bar z\rangle_{\det}$, $\langle \bar a,\bar c : \bar x,\bar z\rangle_{\det}$ and
$\langle \bar a,\bar b : \bar x,\bar y\rangle_{\det}$, so the (triple) determinant product is defined exactly when these three products vanish.

The case $\bar x=\bar a$, $\bar y=\bar b$ makes the connection with commutativity we make in Section \ref{sec: ass vs com in char 2} very explicit. 
One has that
\begin{equation*}
\langle \bar a,\bar b : \bar a,\bar b\rangle_{\det}=\left\{\bar a\,\bar b-(-1)^{|a||b|}\,\bar b\,\bar a\right\}
\end{equation*}
is the graded commutator. 
It always vanishes because $A$ is graded commutative.
\end{remark}

\begin{remark}
\label{rmk: parity}
The parity hypothesis 
\begin{equation*} 
|x|-|a| \;=\; |y|-|b| \;=\; |z|-|c| \equiv 0 \mod 2 
\end{equation*} 
we imposed in Definition \ref{def: determinant Massey product} is what makes the 
determinant of a matrix with two equal columns vanish in the graded setting,
and it cannot be dropped in general. 
Writing $P=|x|-|a|$, one finds
\begin{equation*} 
dU_{\mD}=\bigl(1-(-1)^{P|c|}\bigr)ayc
+\bigl((-1)^{P(|b|+|c|)}-1\bigr)(-1)^{|a||b|}bxc,
\end{equation*} 
so for $P$ odd, 
$U_{\mD}$ is a cocycle when $2=0$ or 
when $|b|$ and $|c|$ are both even, but not in general. 
For instance, for
$|a|=|b|=2$, $|c|=|x|=|y|=3$, $|z|=4$ one computes
\begin{equation*}
dU_{\mD} \;=\;
\begin{vmatrix}
a & x & a \\ b & y & b \\ c & z & c
\end{vmatrix}
\;=\; -2\,acy + 2\,bcx \;\neq\; 0.
\end{equation*}    
\end{remark}

\begin{example}
\label{exa: Iwasaba}
 The running example of this paper will be \emph{Iwasawa manifolds}. 
 These exist at the interface of algebraic topology and differential geometry and are the standard examples
 of non-Kähler manifolds \cite[p. 444]{griffiths78}.
Consider the following groups
\begin{equation*}
  G =  \left\{\begin{pmatrix}
        1 &x &z
        \\
        0 & 1 &y
        \\
        0 &0 &1
    \end{pmatrix} \; \Big| \;  x, y, z \in \Complex 
    \right\}, 
    \qquad  
    \Gamma =  \left\{\begin{pmatrix}
        1 &x &z
        \\
        0 & 1 &y
        \\
        0 &0 &1
    \end{pmatrix} \; \Big| \;  x, y, z \in \Z[i] 
    \right\}.
\end{equation*} 
The group $G$ is the famous Heisenberg group and $\Gamma$ is its integral lattice. 
The group $\Gamma$ acts smoothly on $G$ by left multiplication and thus the  
quotient $\Gamma\backslash G$ forms a compact complex manifold $M$ called the \emph{Iwasawa manifold}. 
Following \cite[Section 2]{Mili24}, the manifold $M$ has a Sullivan minimal model
$A =\Sym( \eta_1, \eta_2, \eta_3, \eta_4, \eta_5, \eta_6)$ with all generators of degree $1$
and differential
\begin{equation*}
d\eta_1=d\eta_2=d\eta_3=d\eta_4=0, \quad d\eta_5 = \eta_1\eta_3- \eta_2\eta_4, \quad d\eta_6 = \eta_2\eta_3 +\eta_1\eta_4.
\end{equation*}
It can be shown that $M$ has vanishing real Massey products of any order with inputs of degree 1,
see \cite[Corollary 2.2]{Mili24}.  
However, the Iwasawa manifold has a nontrivial rational determinant Massey product (Example \ref{exa: Iwasaba formal}).

The cohomology in degree 1 is spanned by the classes $[\eta_1], [\eta_2], [\eta_3], [\eta_4].$
We compute the determinant Massey product 
$\langle [\eta_1], [\eta_2], [\eta_3] \ :\  [\eta_2], - [\eta_1], -[\eta_4] \rangle_{\det}.$
This determinant Massey product has a defining system given 
by choosing our representatives as 
$\left\{ \eta_1, \eta_2, \eta_3, \eta_2, -\eta_1,- \eta_4 \right\}$ and $ \zeta_a 
= -\eta_5, \zeta_b =  \eta_6, \zeta_c = 0$.
We see that 
\begin{equation*}
d\eta_5 = \eta_1\eta_3- \eta_2\eta_4,\quad  d\eta_6 = \eta_2\eta_3 +\eta_1\eta_4
,\quad \textrm{and } \quad d 0 = 0 
=  \eta_1^2 + \eta_2^2.
\end{equation*}
The last formula follows by anticommutativity.
The determinant Massey product associated to this defining system is 
the non-trivial class represented by $U_{\mD}=\eta_1 \eta_5 - \eta_2\eta_6$. 
We finish this computation in Example \ref{exa: Iwasaba formal}.\qed
\end{example}

The next result identifies the indeterminacy subgroup of this product.

\begin{proposition}
\label{prop: indeterminacy subgroup}
The determinant Massey product $\langle \bar a,\bar b,\bar c \ :\ \bar x, \bar y, \bar z\rangle_{\det}$ is well-defined up to the 
indeterminacy subgroup
\begin{equation*}
\bar a\cdot H^{|y|+|c|-1}\left(A\right) + \bar b\cdot H^{|x|+|c|-1}\left(A\right)+ \bar c \cdot H^{|x|+|b|-1}\left(A\right)\;\subseteq\; H^{|a|+|y|+|c|-1}\left(A\right).
\end{equation*}
\end{proposition}

\begin{proof}
The statement asserts that the difference between any two elements in  $\langle \bar a,\bar b,\bar c \ :\ \bar x, \bar y, \bar z\rangle_{\det}$
belongs to the indeterminacy subgroup. 
Let $\overline{U_\mD}$ be a determinant Massey product
represented by a cocycle $U_\mD$ associated to the defining system 
$\mD = \left\{a,b,c,\ \ x,y,z,\ \ \zeta_a,\zeta_b,\zeta_c\right\}$.
Since all determinant Massey products arise from such defining systems,
it suffices to alter the elements of this set and see the relationship between $U_\mD$ and the 
new cocycle $U'$ arising from the altered defining system. 
In Prop. \ref{prop: determinant massey doesnt depend on cocycle representatives},
we proved that a change on the cocycle representatives $a,b,c$ or $x,y,z$ produce cohomologous cocycles.
It remains to check the case where we make a different choice of bounding elements 
$\zeta_a,\zeta_b,\zeta_c$.
This is straightforward: if $\zeta_a'$ is such that $d\zeta_a = d\zeta_a'$, then $\zeta_a-\zeta_a'$ is a 
cocycle; the analogous assertion is true for $\zeta_b-\zeta_b'$ and $\zeta_c-\zeta_c'$.
Therefore, the difference belongs to the indeterminacy subgroup:
\begin{align*}
    U_\mD-U' &= \pm a\left(\zeta_a-\zeta_a'\right) \pm b\left(\zeta_b-\zeta_b'\right) \pm c\left(\zeta_c-\zeta_c'\right).
\end{align*}
\end{proof}

The \emph{determinant Massey products appear in pairs},
because there is an implicit symmetry in their definition.
We explain this fact in the remark below.

\begin{remark} 
\label{rmk: Determinant Massey products appear in pairs}
If the determinant Massey product $\langle \bar a,\bar b,\bar c \ :\ \bar x, \bar y, \bar z\rangle_{\det}$ is defined, 
then we can find representatives 
$a,b,c,x,y,z$ and homogeneous elements $\zeta_a,\zeta_b,\zeta_c$ such that
\begin{equation*}
d\zeta_a=yc-bz, \quad \quad 
d\zeta_b=xc-az, \quad \quad 
d\zeta_c=xb-ay.
\end{equation*}
In this case, we can define the element 
\begin{equation*}
U'_\mD = (-1)^{|x|}x\zeta_a -(-1)^{|y|+|y||x|}y\zeta_b+(-1)^{|z|+|z||x|+|z||y|}z\zeta_c.
\end{equation*}
Assume the underlying algebra is evenly graded, to simplify the signs above.
Then the boundary of the element above is
\begin{align*}
    dU'_{\mD} &= x\left(yc- zb\right) 
    - y\left(xc-za\right) + z\left(xb-ya\right) \\[0.25cm]
    &= xyc - ycx + yza - zya + zxb - xzb=0.
\end{align*}
This way, $U'_{\mD}$ is a cocycle, 
and $-U'_{\mD}$ is the cocycle of the defining system for 
$\langle \bar x,\bar y,\bar z \ : \ \bar a, \bar b, \bar c\rangle_{\det}$ 
formed by the same representatives and $-\zeta_a,-\zeta_b,-\zeta_c$.
Hence $-\overline{U'_{\mD}}$ belongs to $\langle \bar x,\bar y,\bar z \ : \ \bar a, \bar b, \bar c\rangle_{\det}$.
One might be tempted to compare the determinant Massey product sets 
\begin{center}
   $\langle \bar a,\bar b,\bar c \ :\ \bar x, \bar y, \bar z\rangle_{\det}$ \quad and  \quad 
$\langle \bar x,\bar y,\bar z \ : \ \bar a, \bar b, \bar c\rangle_{\det}$. 
\end{center}
However, these sets are not comparable in general. Even when they are both defined, we have: 
\begin{equation*}
|U_\mD| = |a|+|y|+|c|-1, \qquad \textrm{and} \qquad
|U'_\mathcal{D}| = |x|+|y|+|c|-1.
\end{equation*}
To compare these sets, the degrees of $a,b,c$ and $x,y,z$ must match appropriately.
See Proposition \ref{prop: expansion in terms of triple Massey products} for a result precisely in this direction, 
and Proposition~\ref{prop: berglund space},
whose cdga has the two non-trivial and distinct determinant Massey products
$\langle a,b,c : x,y,z\rangle_{\det}=\{[\varepsilon_1]\}$ and 
$\langle x,y,z : a,b,c\rangle_{\det}=\{-[\varepsilon_2]\}$, both of degree $5$.
The reader can check that the boundary of this alternative cocycle is also a determinant, namely
\begin{equation*}
dU_{\mD}' = \begin{vmatrix}
x & x & a \\
y & y & b \\
z & z & c 
\end{vmatrix}.
\end{equation*}
\end{remark}

\begin{example}
In the Iwasawa example, the above pairing is closely related to complexification.  
Interchanging the two triples in the determinant Massey product considered above, 
we see $\left\langle [\eta_2],-[\eta_1],-[\eta_4] \ :\ [\eta_1],[\eta_2],[\eta_3] \right\rangle_{\det}$
has a defining system given by cocycle representatives
$\eta_2,-\eta_1,-\eta_4, \eta_1,\eta_2,\eta_3$
and $(\zeta_a,\zeta_b,\zeta_c)
=(\eta_5,-\eta_6,0).$
The associated determinant Massey product is represented by $U_{\mD'}=-\eta_2\eta_5-\eta_1\eta_6.$

Now, consider the complexification $A\otimes \Complex$,
so that $H^*\left(A \otimes \Complex\right) \cong H^*(A)\otimes \Complex$. 
There, 
\begin{equation*}
d\left(\eta_5+i\eta_6\right) = (\eta_1+i\eta_2)(\eta_3+i\eta_4 ).
\end{equation*}  
The triple Massey product 
$\langle [\eta_1+i\eta_2], [\eta_1+i\eta_2],  [\eta_3+i\eta_4] \rangle$, 
which is non-trivial \cite{Mili24},
contains the class of 
\begin{equation*}
\left(\eta_1 +i \eta_2\right)\left(\eta_5 + i \eta_6\right) = U_\mD + i\left(-U_{\mD'}\right),
\end{equation*}
where $U_\mD = \eta_1\eta_5-\eta_2\eta_6$ is the cocycle of Example \ref{exa: Iwasaba}.
\qed
\end{example}

The next step is to study the behaviour of the determinant Massey products with respect to morphisms of algebras.

\begin{proposition}
\label{prop: properties of the determinant Massey product}
The determinant Massey product $\langle \bar a,\bar b,\bar c \ :\ \bar x, \bar y, \bar z\rangle_{\det}$ is preserved by dg algebra morphisms.
That is, if $\langle \bar a,\bar b,\bar c \ :\ \bar x, \bar y, \bar z\rangle_{\det}$ is defined in $A$
and $f:A\to B$ is a morphism of dg algebras, 
then $\left\langle f^*(\bar a),f^*(\bar b),f^*(\bar c)\ : \ f^*(\bar x), f^*(\bar y), f^*(\bar z)\right\rangle_{\det}$ is defined, and 
\begin{equation*}
f^* \left(\langle \bar a,\bar b,\bar c \ :\ \bar x, \bar y, \bar z\rangle_{\det}\right) \;\subseteq\;
\left\langle f^*(\bar a),f^*(\bar b),f^*(\bar c) \ : \ f^*(\bar x), f^*(\bar y), f^*(\bar z)\right\rangle_{\det}.
\end{equation*}
Furthermore:
\begin{enumerate}
    \item If $f$ is a quasi-isomorphism, then $f^*$ induces a bijection between the determinant Massey products above.
    \item The indeterminacy subgroup is preserved by $f^*$. That is,
    \begin{equation*}
f^*\left(\bar a\cdot H^*(A) + \bar b\cdot H^*(A)+ \bar c \cdot H^*(A)\right) \subseteq 
f^*(\bar a)\cdot H^*(B) + f^*(\bar b)\cdot H^*(B)+ f^*(\bar c) \cdot H^*(B).
\end{equation*}
\end{enumerate}
\end{proposition}

\begin{proof} 
Assume $\langle \bar a,\bar b,\bar c \ :\ \bar x, \bar y, \bar z\rangle_{\det}$ is defined, and fix a 
defining system $\mD=\left\{a,b,c,\ \ x,y,z,\ \ \zeta_a,\zeta_b,\zeta_c\right\}$
for a class $\overline{U_\mD}$ therein.
Let $f:A\to B$ be a morphism of dg algebras.
Since $f$ preserves the algebra structure, it maps the elements in the defining system $\mD$
to elements in $B$ satisfying the corresponding relations in $B$,
giving rise to the element
\begin{align*}
    f(U_\mD) &= f\left((-1)^{|a|}a\zeta_a - (-1)^{|b| + |b||a|}b \zeta_b + (-1)^{|c|+ |c||a| + |c||b| }c\zeta_c\right) \\[0.25cm]
    &= (-1)^{|a|}f(a)f(\zeta_a) - (-1)^{|b| + |b||a|}f(b)f(\zeta_b) + (-1)^{|c|+ |c||a| + |c||b| }f(c)f(\zeta_c).
\end{align*}
Since $f$ commutes with the differentials, the element $f(U_\mD)$ defines a cocycle in $B$.
Furthermore, it represents a determinant Massey product in $\left\langle f^*(\bar a),f^*(\bar b),f^*(\bar c) \ : \ f^*(\bar x), f^*(\bar y), f^*(\bar z)\right\rangle_{\det}$,
which is therefore defined, and moreover $f^*(\overline{U_\mD}) = \overline{f(U_\mD)}$.

That the indeterminacy is preserved is straightforward, it follows from the manipulations done to prove the 
previous item.
It remains to show that if $f$ is a quasi-isomorphism and 
\begin{equation*}
\overline{V}\in \left\langle f^*(\bar a),f^*(\bar b),f^*(\bar c) \ : \ f^*(\bar x),f^*(\bar y), f^*(\bar z)\right\rangle_{\det},
\end{equation*}
then we can find a defining system $\mD$ in $A$ such that the corresponding determinant
Massey product $\overline{U_\mD} \in \left\langle \bar a, \bar b,\bar c \ : \ \bar x,\bar y, \bar z\right\rangle_{\det}$
is defined, and $f^*(\overline{U_\mD}) = \overline{V}$.
Assume $\mathcal E = \left\{a',b',c', \ \ x',y',z', \ \ \eta_a,\eta_b,\eta_c\right\}$ is a defining system for $\overline{V}$,
where without loss of generality, 
\begin{equation*}
a' = f(a), \quad b'=f(b),\quad  \textrm{and}\quad c'=f(c).
\end{equation*}
Our goal is to construct the defining system $\mD$ mentioned before.
Since $x',y',z'$ are cocycles in $B$ and $f^*$ is an isomorphism, there are cocycles $x,y,z$ in $A$ such that 
\begin{equation*}
f(x) = x' +dc_x,\qquad f(y) = y'+dc_y,\quad \textrm{and} \quad f(z)=z'+dc_z.
\end{equation*}
Without loss of generality, we may assume $f(x)=x', f(y)=y',$ and $f(z)=z'$.
Since 
\begin{equation*}
y'c'-b'z' = f\left(yc-bz\right) = d\eta_a,
\end{equation*}
and $b,c,y,z$ are cocycles, there is some $\zeta_a\in A$ with $d\zeta_a = yc-bz$.
Furthermore, since $f^*$ is surjective we can find $\zeta_a$ with the extra property that
\begin{equation*}
f\left(\zeta_a\right) = \eta_a + d\Gamma_a
\end{equation*}
for some $\Gamma_a\in B$.
Similarly, we can find $\zeta_b$ and $\zeta_c$ with $f\left(\zeta_b\right) = \eta_b + d\Gamma_b$ and 
$f\left(\zeta_c\right) = \eta_c + d\Gamma_c$.
The determinant Massey product representative $U_\mD$
associated to the choice $\mD = \left\{a,b,c,\ x,y,z,\ \zeta_a,\zeta_b,\zeta_c \right\}$
satisfies
\begin{align*}
f\left(U_\mD\right) 
& = f\left((-1)^{|a|}a\zeta_a - (-1)^{|b| + |b||a|}b \zeta_b + (-1)^{|c|+ |c||a| + |c||b| }c\zeta_c\right)\\[0.25cm]
&= (-1)^{|a|}a'\left(\eta_a + d\Gamma_a\right) - (-1)^{|b| + |b||a|}b' \left(\eta_b + d\Gamma_b\right) 
+ (-1)^{|c|+ |c||a| + |c||b| }c'\left(\eta_c + d\Gamma_c\right)  \\[0.25cm]
&= V + d\left(a'\Gamma_a - (-1)^{|a||b|}b'\Gamma_b 
+ (-1)^{|b||c| + |a||c| }c'\Gamma_c\right),
\end{align*}
completing the proof.
\end{proof}

A determinant Massey product $\langle \bar a,\bar b,\bar c \ :\ \bar x, \bar y, \bar z\rangle_{\det}$ \emph{vanishes} if it contains the zero class. 
That is, if there is a choice of defining system $\mD$ for it such that the associated cocycle $U_\mD$ is a boundary.

\begin{corollary}
\label{cor: Massey products vanish in formal algebras}
The determinant Massey products vanish in formal commutative dg algebras.
\end{corollary}

\begin{proof}
Recall that a commutative dg algebra $A$ is \emph{formal}
if there is a zig-zag of quasi-isomorphisms of cdgas between $A$ and $H^*(A)$
seen as a cdga with trivial differential.
In view of Prop. \ref{prop: properties of the determinant Massey product},
it suffices to show that the determinant Massey products on a commutative dg algebra with trivial differential
vanish whenever these are defined. 
Indeed, if a determinant Massey product $\langle a,b,c\ : \ x,y,z\rangle_{\det}$
is defined in such an algebra, then we literally have
\begin{equation*}
\langle a,b,c\ : \ x,y,z\rangle_{\det} = a\cdot A^{|y|+|c|-1} + b\cdot A^{|x|+|c|-1} + c\cdot A^{|x|+|b|-1}.
\end{equation*}
In particular, we can choose $\zeta_a=\zeta_b=\zeta_c = 0$, so that the zero element belongs to the set above.
\end{proof}

\begin{example}
\label{exa: Iwasaba formal}
Instantiating to the running Iwasawa example,
the indeterminacy of the determinant Massey product is 
$I = [\eta_1] \cdot H^1(A) +[\eta_2] \cdot H^1(A) +[\eta_3] \cdot H^1(A)$.
Observe that every class in the indeterminacy is represented by an element of $K = \Lambda^2(\eta_1,\eta_2,\eta_3,\eta_4)$.
The subspace of coboundaries $d(A^1) \subset K$ is it is spanned by $d\eta_5 = \eta_1\eta_3- \eta_2\eta_4$ and $d\eta_6 = \eta_2\eta_3 +\eta_1\eta_4$.
As $\eta_1\eta_5-\eta_2\eta_6 \notin K$ the class $[\eta_1\eta_5-\eta_2\eta_6]$ is not in $I$, so the determinant Massey product does not vanish.
We can conclude that the Iwasawa manifold is not formal, and, by the main result of \cite{deligne75} it is also not K\"{a}hler.
\end{example}

We use the following sign convention for the classical triple Massey product. 
If $a,b,c$ are cocycles and $\alpha,\beta$ satisfy $d\alpha=ab$ and $d\beta=bc$,
then $\alpha c-(-1)^{|a|}a\beta$ is the cocycle representative of $\langle \bar a,\bar b,\bar c\rangle$.
With this convention, 
determinant Massey products recover ordinary triple Massey products in the following way.

\begin{proposition} 
\label{prop:ordinary-as-determinant} 
Let $A$ be a cdga and homogeneous elements $a,b,c,x,y,z\in H^*(A)$ be , 
with $ab=0=bc$. 
\begin{itemize}

    \item If $|x|-|a|=|z|-|c|$ is even and $xc=az$, then
    \begin{equation*}
    \langle a,b,c\rangle = -\langle a,0,c : x,b,z \rangle_{\det}. 
    \end{equation*}

    \item If $|y|-|a|=|z|-|c|$ is even and $yc=az$, then
     \begin{equation*}
    \langle a,b,c\rangle = (-1)^{|a||b|}\langle 0,a,c : b,y, z\rangle_{\det}. 
     \end{equation*}
    \item If $|x|-|a|=|y|-|c|$ is even and $xc=ay$, then
     \begin{equation*}
    \langle a,b,c\rangle = (-1)^{|b||c|}\langle a,c,0 : x,y, b\rangle_{\det}. 
    \end{equation*}
\end{itemize}
In particular, the three identities hold when $x,y,z$ are zero classes of the appropriate degrees. 
\end{proposition}

\begin{proof}
This is a direct computation which we leave to the reader.    
\end{proof}

\section{Connections to classical and matric Massey products}
\label{sec: connections to classical and matric Massey products}

In this section, 
we explain some connections of the determinant Massey product with the classical Massey product
and the matric Massey products.

\begin{proposition}
\label{prop: determinant Massey contains Massey}
If a classical triple Massey product $\langle \bar a,\bar b,\bar c\rangle$ is defined on the cohomology of a cdga, 
then the determinant triple Massey product $\langle \bar a,\bar b,\bar c \ : \ \bar a, 0, \bar c\rangle_{\det}$ 
is defined too, 
and 
\begin{equation*}
\langle \bar a,\bar b,\bar c\rangle \subseteq \langle \bar a,\bar b,\bar c \ : \ \bar a, 0, \bar c\rangle_{\det}.
\end{equation*}
\end{proposition}

\begin{proof}
    Assume $\overline m$ is a Massey product in $\langle \bar a,\bar b,\bar c\rangle$ with defining system
    $\{a,b,c,\alpha,\beta\}$.
    That is, 
    \begin{equation*}
    d\alpha = ab, \qquad d\beta = bc, \qquad \textrm{and} \qquad m = \alpha c-(-1)^{|a|}a\beta.
    \end{equation*}
    Then we can take $x=a$, $y=0$, $z=c$, $\zeta_a=-\beta, \zeta_b=0,$ and $\zeta_c=\alpha$. 
    In this case, the bounding conditions of the determinant Massey product become
    \begin{equation*}
    d\zeta_a = -bc, \qquad d\zeta_b=0, \qquad d\zeta_c = ab.
    \end{equation*}
    The determinant Massey product arising from this choice is 
    \begin{align*}
     U_{\mD} &= (-1)^{|a|}a\zeta_a -(-1)^{|b|+|a||b|}b\zeta_b + (-1)^{|c|+|c||a|+|c||b|}c\zeta_c \\[0.25cm]
     &= (-1)^{|a|}a\left(-\beta\right) -0 + (-1)^{|c|+|c||a|+|c||b|}c\alpha \\[0.25cm]
     &= -(-1)^{|a|}a\beta +\alpha c = m.
    \end{align*}
We used that $|\alpha|=|a|+|b|-1$, so that permuting from $c\alpha$ to $\alpha c$ produces the sign that 
precisely cancels the sign in the determinant Massey product.
\end{proof}

Sometimes, the determinant Massey product can be expressed as a linear combination of ordinary triple Massey products.
The next result gives an example of this.

\begin{proposition} 
\label{prop: expansion in terms of triple Massey products} 
Let $\bar a,\bar b,\bar c,\bar x,\bar y,\bar z$ be homogeneous elements in the cohomology of a cdga, 
with $|x|-|a|=|y|-|b|=|z|-|c|$ even, 
such that the triple Massey products $\langle \bar a,\bar y,\bar c\rangle$, 
$\langle \bar a,\bar z,\bar b\rangle$ and $\langle \bar b,\bar x,\bar c\rangle$ are defined.
Then $\langle \bar a,\bar b,\bar c \ : \ \bar x,\bar y,\bar z\rangle_{\det}$ is defined,
and it contains $-u_1+(-1)^{|b||c|}u_2+(-1)^{|a||b|}u_3$ for 
all $u_1\in\langle \bar a,\bar y,\bar c\rangle$, 
$u_2\in\langle \bar a,\bar z,\bar b\rangle$ and $u_3\in\langle \bar b,\bar x,\bar c\rangle$. 
\end{proposition}

\begin{proof} 
Choose representatives $a,b,c,x,y,z$. We compute the three Massey products $\langle \bar a,\bar y,\bar c\rangle$, 
$\langle \bar a,\bar z,\bar b\rangle$ and $\langle \bar b,\bar x,\bar c\rangle$  with defining systems $\{\alpha_1,\beta_1\}$, $\{\alpha_2,\beta_2\}$ and $\{\alpha_3,\beta_3\}$ repectively, i.e., such that $d\alpha_1=ay$, $d\beta_1=yc$, $d\alpha_2=az$, $d\beta_2=zb$, $d\alpha_3=bx$ and $d\beta_3=xc$. 
These have cocycles $U_1=\alpha_1c-(-1)^{|a|}a\beta_1$, $U_2=\alpha_2b-(-1)^{|a|}a\beta_2$ and $U_3=\alpha_3c-(-1)^{|b|}b\beta_3$. 
Then
\begin{equation*}
\zeta_a = \beta_1 - (-1)^{|b||z|}\beta_2, \qquad \zeta_b = \beta_3 - \alpha_2, \qquad \zeta_c = (-1)^{|b||x|} \alpha_3 - \alpha_1 
\end{equation*} 
satisfy the bounding conditions of Definition~\ref{def: determinant Massey product},
and a direct computation, using $|x|\equiv|a|$, $|y|\equiv|b|$ and $|z|\equiv|c|$, 
gives $U_{\mD}=-U_1+(-1)^{|b||c|}U_2+(-1)^{|a||b|}U_3$.
\end{proof}

In \cite{May69}, Peter May generalised ordinary Massey products to \emph{matric 
Massey products} $\langle M_1,M_2,M_3\rangle$, 
whose entries $M_i$ are matrices of cohomology classes.
If $V_1,V_2,V_3$ are matrices with cocycle entries and $A_{12},A_{23}$ are matrices with $dA_{12}=V_1V_2$ and
$dA_{23}=V_2V_3$, then the matric Massey product is represented by the cocycle $A_{12}V_3-\overline{V_1}A_{23}$, 
where $\overline{V_1}$ is obtained from $V_1$ by multiplying each entry $w$ by $(-1)^{|w|}$. 
By the argument of Proposition~\ref{prop: determinant massey doesnt depend on cocycle representatives},
every class of $\langle V_1,V_2,V_3\rangle$ can be computed from any fixed choice of representing matrices.
In the next result, we show that, when $2$ is invertible, a determinant Massey product 
is a matric Massey product.
We revisit the characteristic two case in Section~\ref{sec: ass vs com in char 2}.

\begin{proposition}
\label{prop: connection to matric Massey products}
Suppose that $2$ is invertible in the ground ring. 
Let $\bar a,\bar b,\bar c,\bar x,\bar y,\bar z$ satisfy the degree conditions of 
the determinant Massey products, Definition~\ref{def: determinant Massey product}.
Write $p=|a|$, $q=|b|$, $r=|c|$, $v=(\bar a\ \bar b\ \bar c)$, and let $v^{T}$ be its transpose.
Let 
\begin{equation*} 
W = \begin{pmatrix} 0 & (-1)^{qr}\bar z & -\bar y \\
-(-1)^{pq+pr}\bar z & 0 & (-1)^{pq}\bar x \\ 
(-1)^{pq+pr+qr}\bar y & -(-1)^{pr+qr}\bar x & 0
\end{pmatrix}. 
\end{equation*} 
Then $\langle \bar a,\bar b,\bar c \ : \ \bar x, \bar y, \bar z\rangle_{\det}$ is defined if 
and only if $\langle v, \frac{1}{2}W,v^{T}\rangle$ is defined,
and in that case 
\begin{equation*} 
\langle \bar a,\bar b,\bar c \ : \ \bar x, \bar y, \bar z\rangle_{\det} = \langle v, \tfrac{1}{2}W,v^{T}\rangle. 
\end{equation*} 
\end{proposition}

\begin{proof} 
The proof is a straightforward computation.
First choose representatives $a,b,c,x,y,z$ for the entries of $v$, $\tfrac{1}{2}W$ and $v^T$.
Then the entries of $vW$ and $Wv^T$ are, up to sign, the three minors $yc-bz$, $xc-az$, $xb-ay$.
This gives the first assertion.
Now, suppose that $A_{12}, A_{23}$ are such that 
\begin{equation*}
dA_{12}=\tfrac{1}{2} vW \qquad \textrm{and} \qquad dA_{23}=\tfrac{1}{2} W v^T.
\end{equation*}
Rescaling the entries of $A_{12}$ and $A_{23}$ by $\pm2$ gives triples $\zeta'$ and $\zeta$ satisfying the bounding conditions,
and every such pair arises this way.
An easy computation shows that
\begin{equation*}
A_{12}v^T-\bar vA_{23}=\tfrac12(U_{\zeta'}+U_\zeta)=U_{(\zeta+\zeta')/2}.
\end{equation*}
Since $\tfrac12(\zeta+\zeta')$ also satisfies the bounding conditions, every class in the matric Massey product
is a determinant Massey product. Conversely, taking $\zeta'=\zeta$ shows that every $U_\zeta$ is a matric cocycle.
\end{proof}

\begin{remark} \label{rmk: halving}
The proof of Proposition~\ref{prop: connection to matric Massey products}
uses the invertibility of $2$ only in its last step.
Without this assumption, taking $\zeta' = \zeta$ still gives a defining system for $\langle v,W,v^T\rangle$
with cocycle $2U_\zeta$. 
That is, $U_\zeta$ is literally half of a cocycle representing an element of the matric Massey product $\langle v, W, v^T \rangle$.
In characteristic two this cocycle vanishes, and we show in Section~\ref{sec: ass vs com in char 2}
that $U_\zeta$ need not be a coboundary.
\end{remark}

\section{Applications and examples}
\label{sec: applications}

In this section, we collect three applications of determinant Massey products.

\begin{enumerate}
    \item In Section \ref{sec: ass vs com in char 2}, we explore the behaviour of the determinant Massey product when the ground ring is of characteristic 2.
In this context, determinant Massey products are not matric Massey products 
but a new class of secondary operation.
We give an example of a cdga over $\mathbb F_2$ which is formal as an associative dg algebra, hence with vanishing matric Massey products of all orders, but has a non-trivial determinant Massey product, and so is not formal as a cdga. 
This leads to some pathologies on the level of the deformation complexes.
We conclude by producing a divided power algebra over $\mathbb F_2$ which is formal as an associative dg algebra, but not as a divided power algebra.
This resolves a conjecture of the first author.

\item In Section \ref{sec: Berglund example}, we construct the simplest possible cdga with non-vanishing determinant Massey product and show that it provides an example of a simply-connected finite type rational space which is coformal, with Koszul cohomology ring, vanishing ordinary Massey products of all orders, but which is not formal. 

\item In Section~\ref{sec: Cup-1-algebras},
we extend the construction to cup-1-algebras (Proposition~\ref{prop: cup1prod}). 

\end{enumerate}

\subsection{Associative vs commutative formality in characteristic two}
\label{sec: ass vs com in char 2}

For a graded commutative algebra $H$, denote by $\Hoch^{n,m}(H,H)$ 
the cohomology of the Hochschild cochains with $n$ arguments and internal degree $m$. 
We denote the differential on the Hochschild cochains by $\delta.$
The cochains vanishing on all shuffle products form a subcomplex, the Harrison complex, with cohomology $\Harr^{n,m}(H,H)$. 
The inclusion induces a map $\Harr^{n,m}(H,H) \to \Hoch^{n,m}(H,H)$, which is injective in characteristic zero~\cite{Barr68, Gerst87, ginot10} but in positive characteristic may not be \cite[Section 4]{Barr68}.
We shall show that the ternary operation obtained by homotopy transfer from the determinant Massey product is a Harrison cocycle whose class is a non-zero element of this kernel. 
This element in Hochschild cohomology can also be seen as the so-called \emph{universal Massey product} of $A$ in the sense of \cite{BKS04}.

\begin{proposition}
\label{prop: Harrison not injective}
There is a finite-dimensional cdga $A$ over $\mathbb F_2$, with cohomology $H$, whose universal Massey product
$[m_3] \in \Harr^{3,-1}(H,H)$ is non-zero but maps to zero in $\Hoch^{3,-1}(H,H)$.
This is witnessed by the determinant Massey product in the following sense; there are $a,b,c,x,y,z\in H^2$, such that
\begin{equation*}
\langle a,b,c : x,y,z\rangle_{\det}
=\bigl\{[m_3](a,y,c)+[m_3](a,z,b)+[m_3](b,x,c)\bigr\}\neq\{0\}.
\end{equation*}
In particular, the map
$\Harr^{3,-1}(H,H) \to \Hoch^{3,-1}(H,H)$
induced by the inclusion of complexes is not injective.
\end{proposition}

\begin{proof}
Define a (unital) graded commutative algebra $H$ by
\begin{equation*}
 H^0=\mathbb F_2,\qquad H^2=\mathbb F_2\{a,b,c,x,y,z\},\qquad H^3=\mathbb F_2\{u_1,u_2,u_3\},\qquad H^4=\mathbb F_2\{p_1,p_2,p_3\}, \qquad H^5=\mathbb F_2w,
\end{equation*}
with all other components are zero. Its nonzero products of positive-degree basis elements are
\begin{equation}
\label{ecu: auxiliar}
 bz=cy=p_1,\qquad az=cx=p_2,\qquad ay=bx=p_3,
 \qquad xu_1=yu_2=zu_3=w.
\end{equation}
Associativity follows automatically from the fact that products
of three or more positive-degree elements vanish, as $H^{\geq 6}=0$.

We shall exhibit a nonzero class in the kernel of the inclusion $\Harr^{3,-1} (H,H) \to \Hoch^{3,-1}(H , H)$.
We form a cdga $A$ with cohomology ring $H$ by adjoining two elements $r,s$ of degrees $3,4$ respectively with $\delta r=s$ to $H$.
Set $b z = p_1 + s$ and $a r = w$. 
All other products involving $r$ or $s$ are taken to be zero.
One checks that
\begin{equation*}
 \langle a,b,c:x,y,z\rangle_{\det}=\{[ar]\}=\{w\},  \qquad aH^3+bH^3+cH^3=0.
\end{equation*}
The degree $2$ representatives are unique because any defining system has
$\zeta_a\in r+\mathbb F_2\{u_1,u_2,u_3\}$ and $\zeta_b,\zeta_c\in\mathbb F_2\{u_1,u_2,u_3\}$,
and $aH^3 = bH^3  =cH^3 = \{0\}$.
Transferring along the contraction that has $h(s)=r$ and has zero on all other basis elements, we obtain a ternary operation $m_3$ on $H$. 
This is defined by the property that only nonzero values on the basis are
\begin{equation*}
 m_3(a,b,z) = m_3(a,z,b)=m_3(b,z,a)=m_3(z,b,a)=w.
\end{equation*}
So defined, $m_3$ is a Harrison cocycle. 
To see this, note that since $m_3$ is part of a transferred $A_\infty$-structure with vanishing differential $m_1$, the Stasheff identity for $n=4$ is the condition that $m_3$ is a Hochschild cocycle. 
It is easy to directly check both of the shuffle identities.

This cochain is a coboundary in Hochschild cohomology. 
Explicitly, we have that $\delta \varphi = m_3$, where $\varphi$ is the $2$-cochain of internal degree $-1$ whose only nonzero basis values are
\begin{equation*}
\varphi(a,b) = u_3,\quad \varphi(c,a)=u_2,\quad \varphi(c,b)=u_1,\quad 
\varphi(c,p_3) = \varphi(p_1,a) = \varphi(p_2,b)=w.
\end{equation*}

On the other hand, $[m_3]$ does not vanish in Harrison cohomology.
For a $3$-cochain $f$ of internal degree $-1$, consider $\Phi(f) \coloneqq f(a,y,c ) + f(a,z,b) + f(b,x,c) \in H^5=\mathbb F_2w.$
We shall show that $\Phi$ vanishes on Harrison coboundaries, which implies that it is well defined on Harrison classes.
Let $\psi$ be a Harrison $2$-cochain of internal degree $-1$. 
In characteristic $2$, we have
\begin{equation*}
(\delta \psi)(u,v,t)= u \, \psi(v,t) + \psi(uv,t) + \psi(u,vt) + \psi(u,v) \, t,
\end{equation*}
We have:
\begin{align*}
\Phi(\delta\psi)
&=(\delta\psi)(a,y,c)+(\delta\psi)(a,z,b)+(\delta\psi)(b,x,c)\\
&=a\bigl(\psi(y,c)+\psi(z,b)\bigr)+b\bigl(\psi(x,c)+\psi(a,z)\bigr)+c\bigl(\psi(a,y)+\psi(b,x)\bigr)\\
&\quad+\psi(ay,c)+\psi(a,yc)+\psi(az,b)+\psi(a,zb)+\psi(bx,c)+\psi(b,xc).
\end{align*}
The second line above belongs to $aH^3 + bH^3 + H^3c = 0$, so it vanishes.
The third line also vanishes: in homology, Equation \eqref{ecu: auxiliar} gives  
\begin{equation*}
ay = bx, \qquad cy = bz, \qquad az=cx.
\end{equation*}
Therefore, the third line reduces to $2\psi(ay,c) + 2\psi(cy,bz) + 2\psi(az,cx)$, 
which vanishes because we are in characteristic 2.
We conclude that $\Phi(\delta\psi)=0$.
On the other hand, when $\varphi$ is as above, the same computation gives $\Phi(\delta\varphi)=w$,
consistently with $m_3= \delta \varphi$.
Finally, $\Phi(m_3)=0+w+0=w\neq 0$, which is the determinant Massey product computed above.
\end{proof}

In \cite{campos24}, the authors prove that two commutative dg algebras over a field of characteristic zero are quasi-isomorphic as associative dg algebras if, and only if, they are
quasi-isomorphic as commutative dg algebras. 
The first author showed that this fails in positive characteristic, using higher versions of the Frobenius endomorphism~\cite{Oisin--ObstructionTheory}. 
This left open the possibility that two divided powers algebras over a field of positive characteristic are quasi-isomorphic as associative dg algebras if and only if they are quasi-isomorphic as divided powers algebras; the first author had publicly  conjectured at \emph{Higher Homotopy Algebras in Topology III} in June 2026 that this was true. 
We now show that the cdga $A$ of the proof of Proposition~\ref{prop: Harrison not injective} admits divided powers, which will allow us to use it to negatively resolve this conjecture in characteristic two.
Explicitly, we shall use determinant Massey products to construct an example of a divided powers algebra
that is formal as an associative algebra but not as divided powers algebra.  

\begin{proposition}
\label{prop: dividedpowers}
There is a dg algebra with divided powers over $\mathbb F_2$ that is
formal as an associative dg algebra but not as a commutative dg
algebra. In particular, it is not formal as a dg algebra with divided
powers.
\end{proposition}

\begin{proof}
Let $A$ be the cdga in Proposition~\ref{prop: Harrison not injective}.
For notational convenience, we label the basis numerically as follows: $e_1,e_2,e_3,e_4,e_5,e_6=a,b,c,x,y,z$.
Using the fact that $A$ vanishes in degrees $\geq 6$, to define a system of divided powers on $A$ it suffices to set, for $u=\sum_i\lambda_ie_i\in A^2$, 
\begin{equation*}
\gamma_2(u)=\sum_{i<j}\lambda_i\lambda_je_ie_j,
\end{equation*}
The higher powers $\gamma_k$ are taken to vanish on $A^2$ for $k\geq 3$ and on $A^4$ for $k\geq 2$.
It is easily observed that $\gamma_2(u+v)=\gamma_2(u)+\gamma_2(v)+uv$ and $\gamma_2(\lambda u)=\lambda^2\gamma_2(u)$ as $u^2=0$ for all $u\in A^2$.
Hence this defines a divided powers on $A$.
It is compatible with $d$, because $d$ vanishes on both $A^2$ and on $A^4$.

Recall that $A$ is formal as an associative dg algebra, but not as a commutative dg algebra.
A zig-zag of quasi-isomorphisms of dg algebras with divided powers is in particular one of commutative dg algebras, so $A$ is not formal as a dg algebra with divided powers either.
\end{proof}

\begin{remark}
\label{rmk: associative formality}
The example above shows that the determinant Massey product is not a matric Massey product
in characteristic 2.
Matric Massey products are invariants of the associative homotopy type and so all contain zero on formal associative dg algebras.
\end{remark}

\subsection{An application to Koszul rational spaces}
\label{sec: Berglund example}

We apply determinant Massey products to give an example in rational homotopy theory 
when ordinary Massey products of all orders are not enough.
The main result is Proposition \ref{prop: berglund space}.

\smallskip

The elements $\zeta_a$, $\zeta_b$,
and $\zeta_c$ in a defining system of a determinant Massey product bound quadratic relations.
The fundamental role of quadratic relations in algebra is well established,
see for instance  \cite{Quadraticalgebras}. Among quadratic algebras,
\emph{Koszul algebras} form a subclass which is very well behaved from a homological standpoint.
In \cite{Berglund-KoszulSpaces}, Berglund established an elegant connection of Koszul
algebras to simultaneous formality and coformality in rational homotopy theory.
Recall that a connected, nilpotent, finite $\Q$-type space $X$ is \emph{Koszul}
if it is both formal and coformal.
The main result of \cite{Berglund-KoszulSpaces} asserts the following.

\begin{theorem}[{\cite[Thm. 1.2]{Berglund-KoszulSpaces}}]
Let $X$ be a connected, nilpotent, finite $\Q$-type space. The following are equivalent:
\begin{enumerate}
    \item $X$ is a Koszul space.
    \item $X$ is formal and the cohomology algebra $H^*(X;\Q)$ is Koszul.
    \item $X$ is coformal and the rational homotopy Lie algebra $\pi_*(\Omega X)\otimes \Q$ is Koszul.
\end{enumerate}
\end{theorem}

We use determinant Massey products to show that the result above is sharp in the 
following sense.

\begin{proposition}
\label{prop: berglund space}
There is a simply-connected finite $\Q$-type rational space which is coformal,
with vanishing Massey products of all orders, and Koszul cohomology ring, which is not formal (hence, not Koszul).
An explicit example is the simplicial realization of the quasi-free cdga $A$ on generators $a,b,c,x,y,z$ of degree 2
and $r,s,t$ of degree 3, with differential
\[
 dr=yc-bz,\qquad ds=xc-az,\qquad dt=xb-ay.
\]
Its non-formality
is detected by the determinant Massey product
$\langle a,b,c : x,y,z \rangle_{\det}$.
\end{proposition}

The cohomology of $A$ can be computed explicitly. Consider the degree $5$ cocycles
\[
\varepsilon_1=ar-bs+ct,\qquad \varepsilon_2=xr-ys+zt,
\]
where $\varepsilon_1$ is the representative of $\langle a,b,c : x,y,z \rangle_{\det}$ given by the
defining system $\zeta_a=r$, $\zeta_b=s$, $\zeta_c=t$.
Products of odd-degree classes vanish, so $H^*(A)$ is commutative as a ring, and
there is an isomorphism
\[
H^*(A)\cong\Q[a,b,c,x,y,z,e_1,e_2]/J,
\]
where $J$ is the ideal of the ordinary polynomial ring generated by the $2\times2$ minors of
$\left(\begin{smallmatrix}x&y&z&e_2\\a&b&c&e_1\end{smallmatrix}\right)$, together with
$e_1^2$, $e_1e_2$ and $e_2^2$. One can check that these nine quadrics form a Gr\"obner basis of $J$
for the lexicographic order $x>y>z>e_2>a>b>c>e_1$, so that $H^*(A)$ is Koszul with respect to the
grading in which all the generators have weight $1$ \cite{Fro75,ERT,CDR}.
A similar analysis shows that all Massey products vanish, while
$\langle a,b,c : x,y,z \rangle_{\det}=\{[\varepsilon_1]\}$ is non-trivial. 
The details of these computations will appear elsewhere.

\begin{remark} 
Simply-connected non-formal spaces all of whose Massey products vanish are known.
The point of Proposition~\ref{prop: berglund space} is the combination with coformality and Koszul cohomology.
Since the space $X$ of Prop. \ref{prop: berglund space} 
is coformal but not Koszul, the theorem above shows that $\pi_*(\Omega X)\otimes\Q$ is not a Koszul Lie algebra.
Dually, a formal space that is not coformal and has Koszul homotopy Lie algebra must have non-Koszul cohomology,
because formal spaces with Koszul cohomology are Koszul. 
We do not know whether such a space exists with vanishing higher-order Whitehead products.
It would be interesting to construct one using an analogue of determinant Massey products for dg Lie algebras. 
\end{remark}

\subsection{Cup-1-algebras}
\label{sec: Cup-1-algebras}

Our next application concerns cup-1-algebras.
These algebraic structures have been the subject of recent research, 
as they capture the integral $E_2$-structure of sufficiently nice topological spaces.
In particular, they recover the nilpotent component of the fundamental group for these. 
See \cite[Section 1.1]{Porter--Suciu}, and \cite[Section 4.2.1]{Oisin--ObstructionTheory}.
Extra structure on cochains has been used in a similar way by Porter and Suciu,
who fix part of a defining system with a binomial operation 
to define restricted triple Massey products with smaller indeterminacy 
\cite[Def. 9.8 and Lemma 9.9]{PorterSuciu22}.

\smallskip

A \emph{cup-1-algebra} is a dg algebra $A$ together with a bilinear operation 
$\cup_1\colon A^p\otimes A^q\to A^{p+q-1}$ satisfying the following two identities on homogeneous elements:
\begin{align*} 
\textrm{(Hirsch identity)} 
\quad (ab) \cup_1 c &= a(b \cup_1 c) + (-1)^{|b|(|c| - 1)} \, (a \cup_1 c)b,\\ 
\textrm{(Steenrod identity)} \quad d(a \cup_1 b) &= (da)\cup_1 b +(-1)^{|a|} a\cup_1(db)
+(-1)^{|a|+|b|-1}\bigl(ab -(-1)^{|a||b|}ba\bigr). 
\end{align*} 
A morphism of cup-1-algebras is a morphism of dg algebras that preserves $\cup_1$.

\smallskip

The main result concerning cup-1-algebras appears below.

\begin{proposition}
\label{prop: cup1prod}
Let $A$ be a cup-1-algebra,
and let $a,b,c,x,y,z\in A$ be cocycles with $|a|+|b|=|c|+|x|$, $|b|+|y|=|c|+|z|$ and $|x|+|y|=|a|+|z|$.
Suppose that $\zeta_a,\zeta_b,\zeta_c\in A$ satisfy 
\begin{equation*}
d\zeta_a=ab+cx, \qquad d\zeta_b=by-cz, \qquad d\zeta_c=xy+(-1)^{|a||c|}az.
\end{equation*} 
Then
\begin{equation*}
U = \zeta_ay-(-1)^{|a|}a\zeta_b -(-1)^{|c|} c\zeta_c +(-1)^{|a|+|c|}\left(a\cup_1 c\right)z 
\end{equation*} 
is a cocycle whose cohomology class is well-defined up to the indeterminacy 
$H^{|a|+|b|-1}(A)\cdot\bar y +\bar a\cdot H^{|b|+|y|-1}(A)+\bar c\cdot H^{|x|+|y|-1}(A)$.
Moreover, if $f\colon A\to B$ is a morphism of cup-1-algebras, 
then $f(U)$ is the cocycle obtained in $B$ from $f(a),f(b),f(c),f(x),f(y),f(z)$ 
and $f(\zeta_a),f(\zeta_b),f(\zeta_c)$. 
If $f$ is a quasi-isomorphism, then every class obtained in $B$ from $f(a),f(b),f(c),f(x),f(y),f(z)$
has the form $f^*(\overline{U} )$ for some choice of $\zeta_a,\zeta_b,\zeta_c$ in $A$.
\end{proposition}

\begin{proof}
 Observe that
\begin{align*} 
dU &= (ab+cx)y-a(by-cz)-c\bigl(xy+(-1)^{|a||c|}az\bigr) -\bigl(ac-(-1)^{|a||c|}ca\bigr)z 
\\ 
 &=acz-(-1)^{|a||c|}caz-acz+(-1)^{|a||c|}caz=0. 
\end{align*}
where we have used that the fact that the Steenrod identity implies that
$d(a\cup_1c)=(-1)^{|a|+|c|-1}\bigl(ac-(-1)^{|a||c|}ca\bigr)$.
Any other choice of $\zeta_a$ will differ from $\zeta_a$ by a cocycle $w$, 
which changes $U$ by $wy$, and similarly for $\zeta_b$ and $\zeta_c$.

As $f$ is a morphism of cup-1-algebras, it preserves all operations, 
so $f(\zeta_a),f(\zeta_b),f(\zeta_c)$ satisfy the bounding conditions for the images of the six cocycles,
and the associated cocycle is $f(U)$.
Suppose now that $f$ is a quasi-isomorphism, and let $\zeta_a',\zeta_b',\zeta_c'\in B$ satisfy 
the bounding conditions for $f(a),f(b),f(c),f(x),f(y),f(z)$, with associated cocycle $U'$.
Arguing as in the proof of Proposition~\ref{prop: properties of the determinant Massey product}, 
since $f^*$ is an isomorphism, there are $\zeta_a,\zeta_b,\zeta_c\in A$ satisfying the bounding conditions 
for $a,b,c,x,y,z$, and $\Gamma_a,\Gamma_b,\Gamma_c\in B$, such that
\begin{equation*}
f(\zeta_a)=\zeta_a'+d\Gamma_a, \qquad f(\zeta_b)=\zeta_b'+d\Gamma_b, \qquad f(\zeta_c)=\zeta_c'+d\Gamma_c.
\end{equation*}
If $U$ is the cocycle associated with $\zeta_a,\zeta_b,\zeta_c$, then, 
since the term involving $\cup_1$ is the same in $f(U)$ and in $U'$,
\begin{equation*}
f(U)=U'+d\bigl(\Gamma_a f(y)-f(a)\Gamma_b-f(c)\Gamma_c\bigr).
\end{equation*}
\end{proof}

\begin{remark}
There is a second, related, cocycle, obtained by replacing
$(-1)^{|a|+|c|}\left(a\cup_1 c\right)$ with $(-1)^{|a|+|c|+|a||c|+1}\left(c\cup_1 a\right)$ in the above expression. 
\end{remark}

\begin{proposition}
\label{prop: cup1 and det}
Let $A$ be a cup-1-algebra, and $a,b,c,x,y,z$, $\zeta_a,\zeta_b,\zeta_c$ and $U$
as in Proposition~\ref{prop: cup1prod}. Denote
\[
V_1=\begin{pmatrix} a & c\end{pmatrix},\qquad
V_2=\begin{pmatrix} b & -c\\ x & (-1)^{|a||c|}a\end{pmatrix},\qquad
V_3=\begin{pmatrix} y\\ z\end{pmatrix},
\]
\[
A_{12}=\bigl(\zeta_a,\ (-1)^{|a|+|c|}a\cup_1c\bigr),\qquad
A_{23}=(\zeta_b,\ \zeta_c)^T.
\]
Then:
\begin{enumerate}

\item The element $U$ represents an element of the matric Massey product
$\langle M_1,M_2,M_3\rangle$, where $M_i$ is the matrix of the classes
of the entries of $V_i$.
More precisely, 
the matrices $A_{12}$ and $A_{23}$ satisfy $dA_{12}=V_1V_2$ and
$dA_{23}=V_2V_3$, and $U=A_{12}V_3-\overline{V_1}A_{23}$.

\item Suppose that $A$ is a commutative dg algebra, regarded as a
cup-1-algebra with $\cup_1=0$, and that $|x|-|a|$ is even. Write
$\varepsilon=(-1)^{|a||c|}$. Then the cocycles
$a,c,y,-\varepsilon x,b,z$ and the elements
$\zeta_b,-\varepsilon\zeta_c,-\zeta_a$ form a defining system
$\mD$ of
$\langle \bar a,\bar c,\bar y : -\varepsilon\bar x,\bar b,
\bar z\rangle_{\det}$, and $U_{\mD}=-U$.
\end{enumerate}
\end{proposition}

\begin{proof}
(1) The Steenrod identity for the cocycles $a$ and $c$ gives
$d\bigl((-1)^{|a|+|c|}a\cup_1c\bigr)=(-1)^{|a||c|}ca-ac$. Hence
\[
dA_{12}=\bigl(ab+cx,\ (-1)^{|a||c|}ca-ac\bigr)=V_1V_2,\qquad
dA_{23}=\bigl(by-cz,\ xy+(-1)^{|a||c|}az\bigr)^T=V_2V_3,
\]
and
\[
A_{12}V_3-\overline{V_1}A_{23}
=\zeta_ay+(-1)^{|a|+|c|}(a\cup_1c)z-(-1)^{|a|}a\zeta_b
-(-1)^{|c|}c\zeta_c=U.
\]
The degree conditions of Proposition~\ref{prop: cup1prod} say exactly
that the entries of $V_1V_2$ and $V_2V_3$ are homogeneous.

(2) Both identities defining a cup-1-algebra hold with $\cup_1=0$
because $A$ is graded commutative, so
Proposition~\ref{prop: cup1prod} applies, with
$U=\zeta_ay-(-1)^{|a|}a\zeta_b-(-1)^{|c|}c\zeta_c$. Its degree
conditions give $|b|-|c|=|z|-|y|=|x|-|a|$, which is even, so the
degree conditions of Definition~\ref{def: determinant Massey product}
hold. For this product and the elements
$\zeta_b,-\varepsilon\zeta_c,-\zeta_a$, the bounding conditions read
\[
d\zeta_b=by-cz,\qquad d(-\varepsilon\zeta_c)=-\varepsilon xy-az,\qquad
d(-\zeta_a)=-\varepsilon xc-ab.
\]
The first two are those of Proposition~\ref{prop: cup1prod}, and the
third holds because $xc=(-1)^{|x||c|}cx=\varepsilon cx$, as
$|x|\equiv|a|$. Finally,
$|y|+|y||a|+|y||c|+|y||\zeta_a|\equiv|y|(|b|+|c|)\equiv0$, so
$(-1)^{|y|+|y||a|+|y||c|}y(-\zeta_a)=-\zeta_ay$, and
\[
U_{\mD}=(-1)^{|a|}a\zeta_b+(-1)^{|c|}c\zeta_c-\zeta_ay=-U.
\qedhere
\]
\end{proof}

\section{Higher-arity determinant Massey products}
\label{sec: Higher-order determinant Massey products}

Like their more traditional counterparts, triple determinant Massey products  may be generalised to products of any arity $n\geq 3$.
In this section, we define these products and establish their elementary properties.
We expect most of the structural results of this paper, including the relation to matric Massey products, to have higher-arity analogues.
Some of them, however, are likely to require involved combinatorial arguments, which we leave to future work.

\smallskip

Informally, we start with a $n\times n$ matrix. We choose primitives for the two-by-two minors of the matrix. We use these choices to produce new cocycles; choose primitives for those cocycles and continue the process. At the last stage, keep the cohomology class instead of choosing another primitive to get the value of our operation.

\smallskip

Turning to the explicit computation, let $n\geq 3$ and let $a_1,\ldots,a_n,x_1,\ldots,x_n$ be homogeneous elements of a cdga $A$ such that
\begin{equation}
\label{parity_conditions}
|x_1|-|a_1| \; = \; |x_2|-|a_2| \; = \;\cdots \; = \;|x_n|-|a_n| \; = \;  P \qquad \mbox{ with $P \equiv 0 \pmod{2}$} 
\end{equation}
Consider the square matrix of size $n$ all of whose columns are $(a_1,\ldots,a_n)^T$,
except for the $(n-1)^{th}$ one, which is $(x_1,\ldots,x_n)^T$:
\begin{equation}
\label{ecu: matrix}
M=\begin{pmatrix}
a_1 & \cdots & a_1 & x_1 & a_1 \\
\vdots & & \vdots & \vdots & \vdots \\
a_n & \cdots & a_n & x_n & a_n
\end{pmatrix}.
\end{equation}
We introduce some convenient notation designed to handle the subsets of $[n]:=\{1,\ldots,n\}$
that index the higher-arity version of these products.
For a subset $S\subseteq[n]$ and $k\in[n]$, write
\begin{equation*}
|a_S|=\sum_{t\in S}|a_t|, \quad S_{<k}=\setDef{t\in S}{ t<k} \quad \text{and} \quad
S_{>k}= \setDef{t\in S}{ t >k }.
\end{equation*}
For $k\in S$, denote by $\pos_S(k)\in\{1,\ldots,|S|\}$ the position of $k$
when the elements of $S$ are listed in increasing order.
Suppose $|S| = m\geq 2$, let $M_S$ be the $m \times m$ submatrix of $M$
formed by the rows indexed by $S$ and the last $m$ columns of $M$.
That is to say, $M_S$ is the matrix~\eqref{ecu: matrix} built out the elements with indices in $S$.
If $S=\{s_1<\cdots<s_m\}$, then
\begin{equation*}
M_S=\begin{pmatrix}
a_{s_1} & \cdots & a_{s_1} & x_{s_1} & a_{s_1} \\
\vdots & & \vdots & \vdots & \vdots \\
a_{s_m} & \cdots & a_{s_m} & x_{s_m} & a_{s_m}
\end{pmatrix}.
\end{equation*}
For example, $M_{[n]}=M$.
Clearly, deleting both the first column and the row indexed by $k\in S$ in $M_S$  produces $M_{S\setminus\{k\}}$,
and for $k<l$,
\begin{equation*}
\det\left(M_{\{k,l\}}\right)=\begin{vmatrix} x_k & a_k\\ x_l & a_l\end{vmatrix}= x_ka_l-a_kx_l .
\end{equation*}

Our assumption (\ref{parity_conditions}) immediately implies that the entries in any given row of $M$ all have the same parity.
So, the sign that appears when expanding a minor of $M$ along a column only depends on the row and not on Koszul signs.
We shall use the following formula for the determinant
\begin{equation*}
\sum_{\sigma}\operatorname{sgn}( \sigma)\,m_{ 1\sigma(1)} m_{2\sigma(2)} \cdots m_{m\sigma(m)},
\end{equation*}
where the factors of each term in the above sum is written in the order of the rows.

The first and the last columns of $M_S$ coincide so $\det(M_S)=0$ whenever $|S|\geq 3$.
We recursively define the product via cocycles $U_S$, which are correspond versions to vanishing minors $\det(M_S)$,
just as the cocycle of Definition~\ref{def: determinant Massey product} is a formal version of the vanishing determinant
of a $3\times 3$ matrix with two equal columns.

\begin{definition}
\label{def: higher determinant Massey product}
Let $\bar a_1,\ldots,\bar a_n,\bar x_1,\ldots,\bar x_n\in H^*(A)$ be homogeneous classes
whose degrees satisfy
\begin{equation}
\label{ecu: higher degree condition}
|x_1|-|a_1|  \;=\; |x_2|-|a_2|  \;=\; \cdots  \;=\; |x_n|-|a_n| \;=\; P \qquad \mbox{with } P \equiv 0 \pmod{2}
\end{equation}
A \emph{defining system} for the determinant Massey product of
$\bar a_1,\ldots,\bar a_n$ over $\bar x_1,\ldots,\bar x_n$
is a choice $\mD$ of cocycle representatives $a_1,\ldots,a_n,x_1,\ldots,x_n$ of
$\bar a_1,\ldots,\bar a_n,\bar x_1,\ldots,\bar x_n$, respectively,
together with homogeneous elements $\zeta_S\in A$, one for each $S\subseteq[n]$ with $2\leq |S|\leq n-1$, such that
\begin{equation*}
d\zeta_{\{k,l\}}=\det\left(M_{\{k,l\}}\right)= x_ka_l-a_kx_l \quad \textrm{ for all } \quad 1\leq k<l\leq n,
\end{equation*}
and
\begin{equation*}
d\zeta_S=U_S \quad \textrm{ for all } \quad 3\leq |S|\leq n-1.
\end{equation*}
Here, for every $S\subseteq[n]$ with $|S|\geq 3$,
\begin{equation}
\label{ecu: higher U}
U_S=\sum_{k\in S}(-1)^{\varepsilon(S,k)}\,\zeta_{S\setminus\{k\}}\,a_k,
\end{equation}
with $\varepsilon(S,k)=\pos_S(k)+1+\left|a_k\right|\,\left|a_{S_{>k}}\right|$.
The element $U_{\mD}:=U_{[n]}$ is a cocycle (see Theorem~\ref{thm: higher cocycle}), and we call its class
$\overline{U_{\mD}}$ is the \emph{$n$-fold determinant Massey product} associated to $\mD$.
We denote by
\begin{equation*}
\left\langle \bar a_1,\ldots,\bar a_n\ :\ \bar x_1,\ldots,\bar x_n\right\rangle_{\det}
=\left\{\overline{U_{\mD}} \ \mid\ \mD \text{ is a defining system}\right\}
\subseteq H^{\left |a_{[n]}\right|+P-n+2}(A)
\end{equation*}
the (potentially empty) set of all such classes.
The product is said to be \emph{defined} if this set is non-empty, and to \emph{vanish} if it contains $0$.

\end{definition}
\begin{definition}
    Despite its apparent complexity, the above definition is non-vacuous. One can construct a quasi-free example for each higher products in a similar way as to in Proposition \ref{prop: berglund space}.
    We shall leave the further exploration of properties of such examples to future work.
\end{definition}

We prove the main result of this section next, namely, that the operation produced this way defines a cohomology class.
\begin{theorem}
\label{thm: higher cocycle}
For every defining system $\mD$, the element $U_{\mD}$ is a cocycle.
\end{theorem}

\begin{proof}
We prove that $dU_S=0$ for every $S\subseteq[n]$ with $|S|\geq 3$.
The case $S=[n]$ is precisely the statement.
Since each $a_k$ is a cocycle, the Leibniz rule gives
\begin{equation*}
dU_S=\sum_{k\in S}(-1)^{\varepsilon(S,k)}\,\bigl(d\zeta_{S\setminus\{k\}}\bigr)\,a_k .
\end{equation*}

We have two cases to consider, namely $|S|=3$ and $|S|\geq 4$.
The $|S|=3$ has already been shown immediately after Definition \ref{def: determinant Massey product}.
%
%
When $|S|\ge 4$, each $S\setminus\{i\}$ must have at least three elements. 
 The bounding conditions then give that $d\zeta_{S\setminus\{i\}}=U_{S\setminus\{i\}}$.
In this case, substituting \eqref{ecu: higher U} gives
\begin{equation*}
dU_S=\sum_{i\in S}\ \sum_{j\in S\setminus\{i\}} (-1)^{\varepsilon(S,i)+\varepsilon(S\setminus\{i\},j)}\ \zeta_{S\setminus\{i,j\}}\,a_ja_i .
\end{equation*}
For $i\neq j$, the element $\zeta_{S\setminus\{i,j\}}$ appears exactly twice in the above sum. 
Namely, these are the terms indexed by $(i,j)$ and $(j,i)$.
We show that these two terms carry opposite signs and hence cancel.
Assume without loss of generality that $i<j$, and let $p=\pos_S(i)$ and $q=\pos_S(j)$, so that $p<q$.
Removing $i$ shifts the position of $j$ down by one, but removing $j$ does not change the position of $i$, i.e., $\pos_{S\setminus\{i\}}(j)=q-1$ and $\pos_{S\setminus\{j\}}(i)=p.$
Similarly, we have $(S\setminus\{i\})_{>j}=S_{>j}$ and $(S\setminus\{j\})_{>i}=S_{>i}\setminus\{j\}$.
Writing $\eta=|a_i||a_{S_{>i}}|+|a_j||a_{S_{>j}}|$, the definition of $\varepsilon$ gives
\begin{equation*}
\varepsilon(S,i)+\varepsilon(S\setminus\{i\},j)\equiv p+q+1+\eta,
\qquad\quad
\varepsilon(S,j)+\varepsilon(S\setminus\{j\},i)\equiv p+q+\eta+|a_i||a_j| \pmod 2.
\end{equation*}
Koszul signs give us that $a_ia_j=(-1)^{|a_i||a_j|}a_ja_i$, and therefore the two terms are
\begin{equation*}
(-1)^{p+q+1+\eta}\,\zeta_{S\setminus\{i,j\}}\,a_ja_i
\qquad\text{and}\qquad
(-1)^{p+q+\eta+|a_i||a_j|}\,\zeta_{S\setminus\{i,j\}}\,a_ia_j=(-1)^{p+q+\eta}\,\zeta_{S\setminus\{i,j\}}\,a_ja_i,
\end{equation*}
which cancel as they have opposite signs.
\end{proof}


The following result collects together the basic properties of this product. For $S= \{s_1<\cdots<s_m\}\subseteq[n]$, we write
$\langle \bar a_{S}:\bar x_{S}\rangle_{\det}$ to denote the product $\langle \bar a_{s_1},\ldots,\bar a_{s_m} 
: \bar x_{s_1},\ldots,\bar x_{s_m}\rangle_{\det}$.

\begin{proposition}
\label{prop: higher properties}
Let $\bar a_1,\ldots,\bar a_n,\bar x_1,\ldots,\bar x_n$ be as in
Definition~\ref{def: higher determinant Massey product},
and assume that the determinant Massey product 
$\langle \bar a_1, \dots,\bar a_n:\bar x_1,\dots,\bar x_n \rangle_{\det}$ is defined.
The following holds:
\begin{enumerate}
\item 
Every $\langle \bar a_{S}:\bar x_{S}\rangle_{\det}$ with $S\subsetneq[n]$ and $|S|\geq 3$ is defined and vanishes.
More precisely, they vanish coherently, a defining system $\mD$ for
$\langle \bar a_1,\ldots,\bar a_n:\bar x_1,\ldots,\bar x_n\rangle_{\det}$
restricts to a defining system for each $\langle \bar a_{S}:\bar x_{S}\rangle_{\det}$ whose cocycle $U_S$ is a coboundary.
\item 
If $f\colon A\to B$ is a morphism of cdgas, then $\langle f^*(\bar a_1),\ldots,f^*(\bar a_n) : f^*(\bar x_1),\ldots,f^*(\bar x_n)\rangle_{\det}$ is defined and
\begin{equation*}
f^*\langle \bar a_1,\ldots,\bar a_n:\bar x_1,\ldots,\bar x_n\rangle_{\det}\subseteq
\langle f^*(\bar a_1),\ldots,f^*(\bar a_n) : f^*(\bar x_1),\ldots,f^*(\bar x_n)\rangle_{\det}.
\end{equation*}
If $f$ is a quasi-isomorphism, the inclusion is a bijection.
\item 
The product set is a union of cosets to the subgroup
$\sum_{k=1}^n H^{|a_{[n]}|-|a_k|+ P -n+2}(A)\cdot \bar a_k$.
\item 
If $A$ has zero differential, the product vanishes.
In particular, this occurs if $A$ is formal.
\end{enumerate}
\end{proposition}

\begin{proof}
Let $\mD$ be a defining system.
\begin{enumerate}
\item 
This is straightforward to check from the definitions.

\item 
The morphism $f$ maps defining systems to defining systems, and $f(U_{\mD})=U_{f(\mD)}$.
If $f$ is a quasi-isomorphism, then the desired bijection is constructed by induction on $|S|$,
via generalising version of the proof of Proposition~\ref{prop: properties of the determinant Massey product}.

\item
Replacing $\zeta_{[n]\setminus\{k\}}$ by $\zeta_{[n]\setminus\{k\}}+w$, with $w$ a cocycle,
gives a defining system whose cocycle is $U_{\mD}+(-1)^{\varepsilon([n],k)}\,w\,a_k$.

\item 
If $d=0$, the bounding conditions force $x_ka_l=a_kx_l$, and then the choice $\zeta_S=0$ for all $S$
is a defining system with $U_{\mD}=0$.
If $A$ is formal, the claim follows by invariance under quasi-isomorphisms,
applied along a zig-zag connecting $A$ and $H^*(A)$.
\end{enumerate}
\end{proof}

We expect these higher-arity products to be related to matric Massey products when $2$ is invertible in the ground ring.

\smallskip

\noindent {\bf{Acknowledgments:}}
The first author received funding from the European Union’s Horizon 2020 research and innovation programme under the Marie Skłodowska-Curie grant agreement No 945322. 
He also appreciates the hospitality of IMUS, Seville and Leiden University.
The second author has been partially supported by the MICINN grant PID2023-149804NB-I00,
the Junta de Andalucía grant ProyExcel-00827, 
and the Junta de Andalucía research group grupoPaidi-G-FEDER-FQM264. 
We thank A. Berglund for useful conversations. 

\medskip

\noindent {\bf{AI use disclosure:}}
This document has been in preparation since early 2025, and the central idea originated during a February 2024 research trip.
The authors completed the first draft independently of generative AI. 
LLMs were subsequently used to improve the presentation of the text and to perform specific 
computations involved in preparing parts of the examples in Sec. \ref{sec: applications}.

\bibliographystyle{plain}
\bibliography{MyBib}

\bigskip

\noindent\sc{Oisín Flynn-Connolly}\\ 
\noindent\sc{Leiden Institute of Advanced Computer Science \\ Leiden University, Leiden, The Netherlands }\\
\noindent\tt{oisinflynnconnolly@gmail.com}\\

\noindent\sc{José Manuel Moreno Fernández}\\ 
\noindent\sc{Departamento de Álgebra, Geometría y Topología \\ Universidad de Málaga, 29080, Málaga, Spain}\\
\noindent\tt{josemoreno@uma.es}\\

\noindent\sc{Fernando Muro}\\
\noindent\sc{Departamento de Álgebra, Universidad de Sevilla \\ 41012 Sevilla, Spain}\\
\noindent\tt{fmuro@us.es}\\
 \end{document}